\documentclass{amsart}

\usepackage{a4wide,latexsym,amssymb,amsthm,amsmath}

\theoremstyle{plain}

\newtheorem{theorem}{Theorem}
\newtheorem{lemma}{Lemma}
\newtheorem{proposition}{Proposition}

\theoremstyle{remark}

\newtheorem{remark}{Remark}
\newtheorem{example}{Example}

\def \Rt {\widetilde{R}}

\def\R{\mathbb{R}}
\def\C{\mathbb{C}}
\def\E{\mathbb{E}}
\def\H{\mathbb{H}}
\def\S{\mathbb{S}}

\newcommand{\Sol}{\mathrm{Sol}_3}

\begin{document}

\title[Gauss maps of constant mean curvature surfaces]
{Gauss maps of constant mean curvature surfaces in three-dimensional
homogeneous spaces}

\author[J. Inoguchi]{Jun-ichi Inoguchi}
\address{Yamagata University\\ Department of Mathematical Sciences,
Faculty of Science\\ Kojirakawa-machi 1-4-12\\ Yamagata 990-8560\\
Japan} \email{inoguchi@sci.kj.yamagata-u.ac.jp}

\author[J. Van der Veken]{Joeri Van der Veken}
\address{Katholieke Universiteit Leuven\\ Departement
Wiskunde\\ Celestijnenlaan 200 B -- Box 2400\\ BE-3001 Leuven\\
Belgium} \email{joeri.vanderveken@wis.kuleuven.be}
\thanks{The second author is a postdoctoral researcher supported by the Research Foundation -- Flanders (F.W.O.)}

\begin{abstract}
It is well-known that for a surface in a 3-dimensional real space
form the constancy of the mean curvature is equivalent to the
harmonicity of the Gauss map. However, this is not true in general
for surfaces in an arbitrary 3-dimensional ambient space. In this
paper we study this problem for surfaces in an important and very
natural family of 3-dimensional ambient spaces, namely homogeneous
spaces. In particular, we obtain a full classification of constant
mean curvature surfaces, whose Gauss map satisfies the more mild
condition of vertical harmonicity, in \emph{all} 3-dimensional
homogeneous spaces.
\end{abstract}



\maketitle

\section{Introduction}
In recent years, there has been done a lot of research concerning
surfaces in 3-dimensional homogeneous spaces. Initial work was done
in \cite{IKOS}. In particular, the study of constant mean curvature
surfaces (CMC surfaces) and minimal surfaces in 3-dimensional
homogeneous spaces is paid much attention to by differential
geometers. We mention for example the paper \cite{Daniel}.

One of the reasons for this success is that homogeneous spaces are
among the most natural candidates for the role of ambient space in
submanifold theory and that in dimension 3 the classification of
these spaces is well-understood. Indeed, we have the following.
\begin{theorem} \label{theo1} Let $M^3$ be a 3-dimensional simply connected homogeneous Riemannian manifold with
isometry group $I(M^3)$, i.e. $I(M^3)$ acts transitively on $M^3$.
Then $\dim I(M^3)\in\{3,4,6\}$ and moreover:
\begin{itemize}
\item[(i)] if $\dim I(M^3)=6$, then $M^3$ is a real space form of constant sectional curvature $c$,
i.e. Euclidean space $\E^3$, a hyperbolic space $\mathbb{H}^3(c)$ or
a three-sphere $\mathbb{S}^3(c)$,
\item[(ii)] if $\dim I(M^3)=4$, then $M^3$ is locally isometric to a Bianchi-Cartan-Vranceanu space (different from
$\E^3$ and $\mathbb{S}^3(c)\setminus\{\infty\}$), i.e. $M^3$ is a
Riemannian product $\mathbb{H}^2(c)\times\R$ or
$\mathbb{S}^2(c)\times\R$, the Heisenberg group equipped with a
left-invariant metric $\mathrm{Nil}_3$, or one of following Lie
groups equipped with a left-invariant metric yielding a
four-dimensional isometry group: the special unitary group
$\mathrm{SU}(2)$ or the universal covering of the special linear
group $\widetilde{\mathrm{SL}}(2,\R)$,
\item[(iii)] if $\dim I(M^3)=3$, then $M^3$ is a general
3-dimensional Lie group with left-invariant metric.
\end{itemize}
\end{theorem}
Another reason is that the classification above contains the eight
model geometries of Thurston (cfr. \cite{Thurston}), namely $\E^3$,
$\H^3$, $\S^3$, $\H^2\times\R$, $\S^2\times\R$,
$\widetilde{\mathrm{SL}}(2,\mathbb{R})$, $\mathrm{Nil}_3$ and
$\Sol$. The famous geometrization conjecture of Thurston states that
these eight spaces are the `building blocks' to construct \emph{any}
3-dimensional manifold.

Gauss maps of surfaces in 3-dimensional real space forms play an
important role in surface geometry. In fact, as shown by E.~Ruh and
J.~Vilms in \cite{RV}, for submanifolds of a Euclidean space the
harmonicity of the Gauss map is equivalent to the fact that the mean
curvature vector is parallel in the normal bundle. For surfaces in
$\E^3$, this corresponds to the constancy of the mean curvature. In
particular, the minimality of a surface in $\E^3$ is equivalent to
the holomorphicity of Gauss map. The characterization due to
Ruh-Vilms was generalized to CMC surfaces in 3-dimensional real
space forms (i.e. in spaces of type (i) in Theorem \ref{theo1}) by
T.~Ishihara in \cite{Ishihara}.

On the contrary, in 3-dimensional homogeneous spaces of
non-constant curvature, the harmonicity of Gauss map is a very
strong restriction for CMC surfaces. In fact, A. Sanini (cfr.
\cite{Sanini}) and M. Tamura (cfr. \cite{Tamura}) showed that the
only CMC surfaces with harmonic Gauss map in a 3-dimensional
homogeneous space with 4-dimensional isometry group are inverse
images of geodesics under the Hopf-fibration or totally geodesic
leaves. The latter case only occurs if the ambient space is a
direct product space. Instead of harmonicity, Sanini and Tamura
studied a more mild condition, namely ``vertical harmonicity" of
the Gauss map for CMC surfaces in 3-dimensional homogeneous spaces
with 4-dimensional isometry group (i.e. in spaces of type (ii) of
Theorem \ref{theo1}).

In this paper we classify CMC surfaces with vertically harmonic
Gauss map in 3-dimensional homogeneous spaces with 3-dimensional
isometry group (i.e. in spaces of type (iii) of Theorem
\ref{theo1}). For this purpose, we have to treat \emph{all}
3-dimensional Lie groups, for which we use methods based on results
of J.~Milnor from \cite{Milnor}. Our result finishes the
classification of CMC surfaces with vertically harmonic Gauss map in
\emph{all} 3-dimensional homogeneous spaces.

\section{CMC surfaces and Gauss maps}

\subsection{CMC hypersurfaces}
Consider Riemannian manifolds $(M^n,g_M)$ and $(Q^{n+1},g_Q)$ with
Levi Civita connections $\nabla^M$ and $\nabla^Q$ respectively. Let
$f:M^n\rightarrow Q^{n+1}$ be an isometric immersion with unit
normal $N$. The \emph{second fundamental form} $h$ is a field of
symmetric bilinear forms on $M^n$, defined by
$$\nabla_{\mathrm{d}f(X)}^Q \mathrm{d}f(Y)=\mathrm{d}f(\nabla^M_XY)+h(X,Y)N$$
for vector fields $X,Y$ on $M^n$. The \emph{shape operator} $S$ is a
field of symmetric operators on $M^n$, defined by
$SX=-\mathrm{d}f^{-1}(\nabla^Q_{\mathrm{d}f(X)}N)$ for $p\in M^n$
and $X\in T_pM^n$, and it is related to the second fundamental form
by $h(X,Y)=g_M(SX,Y)=g_M(X,SY)$ for $p\in M^n$ and $X,Y\in T_pM^n$.
The mean curvature $H$ is a function on $M^n$ defined by
$$H(p)=\mathrm{tr}_{g_M|_p}(h|_p).$$ The immersion is said to have
\emph{constant mean curvature} (CMC) if $H$ is constant, and it is
said to be \emph{minimal} if $H$ vanishes identically.

\subsection{Grassmannian bundle}
Let $(Q^n,g_Q)$ be a Riemannian $n$-manifold. Denote by
$\mathrm{Gr}_{\ell}(TQ^n)$ the \textit{Grassmannian bundle} of
$\ell$-planes in the tangent bundle $TQ^n$:
$$
\mathrm{Gr}_{\ell}(TQ^n):= \bigcup_{q\in Q^n}
\mathrm{Gr}_{\ell}(T_{q}Q^n).
$$
Here $\mathrm{Gr}_{\ell}(T_{q}Q^n)$ denotes the \textit{Grassmannian
manifold} of $\ell$-planes in the tangent space
$(T_{q}Q^n,g_{Q}\vert_{q})$ at $q$. The Grassmannian bundle
$\mathrm{Gr}_{\ell}(TQ^n)$ is a fiber bundle over $Q^n$ associated
to the orthonormal frame bundle $\mathrm{O}(Q^n)$. The standard
fiber of this bundle is the Grassmannian manifold
$\mathrm{Gr}_{\ell}(\mathbb{R}^n)$. The canonical $1$-form of
$\mathrm{O}(Q^n)$ and the Levi Civita connection of $g_Q$ induce an
invariant Riemannian metric on $\mathrm{Gr}_{\ell}(TQ^n)$. With
respect to this metric, the natural projection
$\pi:\mathrm{Gr}_{\ell}(TQ^n)\to  Q^n$ is a Riemannian submersion
with totally geodesic fibers. For more details on the Riemannian
structure of $\mathrm{Gr}_{\ell}(TQ^n)$, we refer to \cite{JR} and
\cite{Sanini}.

Let $f:M^m\to Q^n$ be an immersed submanifold. Then the
\textit{(tangential) Gauss map} $\psi$ of $f$ is a smooth map of
$M^m$ into $\mathrm{Gr}_{m}(TQ^n)$ defined by
$$\psi(p):=\mathrm{d}f_{p}(T_{p}M^m)\in \mathrm{Gr}_{m}(T_{f(p)}Q^n).$$

\begin{remark}
The case $(Q^3,g_Q)=\mathbb{E}^3$ is exceptional. In fact, since
$\mathbb{R}^3$ has absolute parallelism, the Grassmannian bundle is
a trivial fiber bundle:
$\mathrm{Gr}_{2}(T\mathbb{R}^3)=\mathbb{R}^{3}\times \mathrm{Gr}_{2}
(\mathbb{R}^3)=\mathbb{R}^{3}\times \mathbb{R}P^2$. If we consider
the Grassmannian bundle $\mathrm{Gr}_{2}^{+}(T\mathbb{R}^3)$ of all
\emph{oriented} $2$-planes, we have
$\mathrm{Gr}^{+}_{2}(T\mathbb{R}^3)=\mathbb{R}^{3}\times
\mathbb{S}^2$. If $f:M^2\rightarrow\E^3$ is an isometric immersion
of a surface with unit normal $N$, then the oriented tangential
Gauss map $\psi:M^2\to \mathrm{Gr}_{2}^{+}(T\mathbb{R}^3)$ is given
by $\psi=(f,N)$. Hence we may ignore the first component and we
obtain the classical Gauss map $\psi=N:M^2\rightarrow\S^2$.
\end{remark}

\begin{remark}
For a non-empty subset $\Sigma$ of $\mathrm{Gr}_{\ell}(TQ^n)$, the
totality of $\ell$-submanifolds all of whose tangent spaces belong
to $\Sigma$ is called the $\Sigma$-geometry. The \textit{Grassmann
geometry} on $Q^n$ is the collection of such $\Sigma$-geometries in
$Q^n$. Let $G$ be the identity component of the isometry group of
$Q^n$. Then $G$ acts isometrically on $\mathrm{Gr}_{\ell}(TQ^n)$. If
$\Sigma$ is a $G$-orbit in $\mathrm{Gr}_{\ell}(TQ^n)$, then the
$\Sigma$-geometry  is said to be of \textit{orbit type}. Note that
if $Q^n$ is a homogeneous Riemannian manifold, then an orbit
$\Sigma$ is a subbundle of $\mathrm{Gr}_{\ell}(TQ^n)$. H.~Naitoh has
developed the general theory of Grassmann geometries of orbit type
on Riemannian symmetric spaces. Grassmann geometries on the
3-dimensional Heisenberg group and on the motion groups
$\mathrm{E}(1,1)$ and $\mathrm{E}(2)$ are investigated by Naitoh,
Kuwabara and the first named author in \cite{IKN} and
\cite{Kuwabara}.
\end{remark}

\subsection{Harmonic maps}
Next, we recall some fundamental ingredients of harmonic map theory
from the lecture notes \cite{EL}.

Let $(M^m,g_M)$ and $(Q^n,g_Q)$ be Riemannian manifolds, with Levi
Civita connections $\nabla^M$ and $\nabla^Q$ respectively. Let
$f:M^m\to Q^n$ be a smooth map. Then the \textit{energy density}
$e(f)$ of $f$ is a function on $M^m$ defined by
$e(f)=|\mathrm{d}f|^{2}/2$. One can see that $f$ is constant if and
only if $e(f)=0$ and that $e(f)=m/2$ if $f$ is an isometric
immersion. The \textit{energy} $E(f;\mathfrak{D})$ of $f$ over a
region $\mathfrak{D}\subset M^m$ is
$$E(f;\mathfrak{D})=\int_{\mathfrak D}e(f)\>\mathrm{d}v_{M}.$$
A smooth map $f$ is said to be \textit{harmonic} if it is a critical point
of the energy over every compactly supported region of $M$.

The second fundamental form $\nabla \mathrm{d}f$ of $f$ is in fact
an extension of the second fundamental form for isometric immersions
and is defined by
$$\nabla \mathrm{d}f(X,Y)=
\nabla^{Q}_{\mathrm{d}f(X)}\mathrm{d}f(Y)-\mathrm{d}f(\nabla_{X}^{M}Y)$$
for vector fields $X,Y$ on $M^m$. The second fundamental form
$\nabla \mathrm{d}f$ is a symmetric $TQ^n$-valued tensor field,
i.e., $(\nabla \mathrm{d}f)(X,Y)=(\nabla \mathrm{d}f)(Y,X)$. The
trace $\tau(f):=\mathrm{tr}_{g_M}(\nabla \mathrm{d}f)$ is called the
\textit{tension field} of $f$. It is known that $f$ is harmonic if
and only if $\tau(f)=0$. Remark that an isometric immersion is
harmonic if and only if it is minimal.

\subsection{Vertically harmonic maps}
Let $(P,g_P)$ be a Riemannian manifold and $\pi:(P,g_P)\to (Q,g_Q)$
a Riemannian submersion. With respect to the metric $g_P$, the
tangent bundle $TP$ of $P$ has a splitting
\begin{equation}\label{Splitting}
T_{u}P=\mathcal{H}_{u}\oplus \mathcal{V}_{u}, \ \ u\in P.
\end{equation}
Here $\mathcal{V}_{u}=\mathrm{Ker}\>(\mathrm{d}\pi_{u})$ and
$\mathcal{H}_{u}$ is the orthogonal complement of $\mathcal{V}_{u}$
in $T_{u}P$. The linear subspaces $\mathcal{V}_{u}$ and
$\mathcal{H}_{u}$ are called the \textit{vertical subspace} and the
\textit{horizontal subspace} of $T_{u}P$.

Now let $f:M\to P$ be a smooth map. Then its tension field $\tau(f)$ is
decomposed as
$$\tau(f) = \tau^{\mathcal H}(f)+\tau^{\mathcal V}(f)$$
according to the splitting (\ref{Splitting}). A smooth map $f$ is
said to be \textit{vertically harmonic} if the vertical component
$\tau^{\mathcal V}(f)$ of the tension field vanishes. C.~M.~Wood has
shown that vertical harmonicity of $f$ is equivalent to the
criticality of $f$ with respect to the vertical energy through
vertical variations in \cite{cW}.

As announced in the introduction, the purpose of the present paper
is to study vertical harmonicity of Gauss maps for CMC surfaces in
3-dimensional Lie groups with left-invariant metric, yielding a
3-dimensional isometry group. For our purpose, we recall the
following result due to Sanini.
\begin{lemma}{\rm(\cite{Sanini})} \label{lem3.1}
Let $(Q^3,g_Q)$ be a Riemannian $3$-manifold and $f:M^2\rightarrow
Q^3$ a CMC surface with unit normal $N$. Take a principal frame
field
$\{\epsilon_1=\mathrm{d}f(E_1),\epsilon_2=\mathrm{d}f(E_2),\epsilon_3=N\}$,
where $\{E_1,E_2\}$ is an orthonormal frame field on $M^2$ which
diagonalizes the shape operator associated to $N$. Put
$\mathfrak{R}_{ijkl}:=g_{Q}(R^Q(\epsilon_i,\epsilon_j)\epsilon_k,\epsilon_l).$
Then we have the following.
\begin{itemize}
\item[(i)] The Gauss map $\psi$ is vertically harmonic if and only if
$\mathfrak{R}_{1213}=\mathfrak{R}_{2123}=0$. Moreover, when $f$ is
minimal, $\psi$ is harmonic if and only if, in addition,
$\mathfrak{R}_{3113}=\mathfrak{R}_{3223}$.
\item[(ii)] $\psi$ is conformal if and only if $f$ is minimal or totally umbilical.
\end{itemize}
\end{lemma}

The following statement is an immediate corollary of Lemma
\ref{lem3.1}.
\begin{lemma}\label{lem3.2}
Let $(Q^3,g_Q)$ be a Riemannian $3$-manifold and let
$f:M^2\rightarrow Q^3$ be a CMC surface. The Gauss map is vertically
harmonic if and only if the normal component of
$R^Q(\mathrm{d}f(X),\mathrm{d}f(Y))\mathrm{d}f(Z)$ vanishes for all
$p\in M^2$ and for all $X,Y,Z\in T_pM^2$.
\end{lemma}

In \cite{IVdV2}, we have proven the following lemma.
\begin{lemma}{\rm(\cite{IVdV2})}\label{lem3.3}
Consider a Riemannian manifold $(Q^3,g_Q)$ and an orthonormal frame
field $\{e_1,e_2,e_3\}$ on $Q^3$ which diagonalizes the Ricci
tensor. Denote by $K_{ij}$ the sectional curvature of the plane
spanned by $\{e_i,e_j\}$, for $i,j\in\{1,2,3\}$. Let
$f:M^2\rightarrow Q^3$ be a surface such that the normal component
of $R^Q(\mathrm{d}f(X),\mathrm{d}f(Y))\mathrm{d}f(Z)$ vanishes
identically for all vector fields $X,Y,Z$ on $M^2$ and suppose that
$N=\alpha e_1+\beta e_2+\gamma e_3$ is a unit normal. Then every
point of $M^2$ has an open neighbourhood in $M^2$ on which at least
on the of the following holds:
\begin{itemize}
\item[(i)] $\alpha=\beta=0$,
\item[(ii)] $\alpha=\gamma=0$,
\item[(iii)] $\beta=\gamma=0$,
\item[(iv)] $\alpha=0$ and $K_{12}=K_{13}$,
\item[(v)] $\beta=0$ and $K_{12}=K_{23}$,
\item[(vi)] $\gamma=0$ and $K_{13}=K_{23}$,
\item[(vii)] $K_{12}=K_{23}=K_{13}$.
\end{itemize}
\end{lemma}

\section{Three-dimensional Lie groups}\label{sec3}

\subsection{Unimodular Lie groups}
A Lie group $G$ is said to be \textit{unimodular} if its
left-invariant Haar measure is right-invariant. We refer to the work
of Milnor \cite{Milnor} for an infinitesimal reformulation of the
unimodularity property.

Given a 3-dimensional unimodular Lie group $G$ with a left-invariant
metric, there exists an orthonormal basis $\{e_1,e_2,e_3\}$ for the
Lie algebra $\mathfrak{g}$, satisfying
\begin{equation}\label{Lieunimod}
[e_1,e_2]=c_{3}e_{3},\quad  [e_2,e_3]=c_{1}e_{1},\quad
[e_3,e_1]=c_{2}e_{2}, \qquad c_{i}\in \mathbb{R}.
\end{equation}
To describe the Levi Civita connection and the curvature of $G$, we
introduce the following constants:
$$
\mu_{i}=\frac{1}{2}(c_{1}+c_{2}+c_{3})-c_{i}.
$$
\begin{proposition}\label{prop2}
Let $G$ be a 3-dimensional unimodular Lie group with a
left-invariant metric $\langle \cdot,\cdot\rangle$ and use the
notations introduced above. The Levi Civita connection
$\widetilde{\nabla}$ of $G$ is given by
$$
\begin{array}{lll}
\widetilde{\nabla}_{e_1}e_{1}=0, & \widetilde{\nabla}_{e_1}e_{2}=\mu_1e_{3}, & \widetilde{\nabla}_{e_1}e_{3}=-\mu_{1}e_{2},\\
\widetilde{\nabla}_{e_2}e_{1}=-\mu_{2}e_{3}, & \widetilde{\nabla}_{e_2}e_{2}=0, & \widetilde{\nabla}_{e_2}e_{3}=\mu_{2}e_{1},\\
\widetilde{\nabla}_{e_3}e_{1}=\mu_{3}e_2, &
\widetilde{\nabla}_{e_3}e_{2}=-\mu_{3}e_{1}, &
\widetilde{\nabla}_{e_3}e_{3}=0.
\end{array}
$$
The frame $\{e_1,e_2,e_3\}$ diagonalizes the Ricci tensor and the
Riemann-Christoffel curvature tensor $\widetilde{R}$ is determined
by the following sectional curvatures:
\begin{eqnarray*}
&& K_{12} = \langle \widetilde{R}(e_1,e_2)e_2,e_1\rangle = c_{3}\mu_{3}-\mu_1\mu_2,\\
&& K_{23} = \langle \widetilde{R}(e_2,e_3)e_3,e_2\rangle = c_{1}\mu_{1}-\mu_{2}\mu_{3},\\
&& K_{13} = \langle \widetilde{R}(e_1,e_3)e_3,e_1\rangle =
c_{2}\mu_{2}-\mu_{1}\mu_{3}.
\end{eqnarray*}
\end{proposition}

Milnor classified 3-dimensional unimodular Lie groups based on the
signature of the constants $(c_1,c_2,c_3)$.
\begin{center}
\begin{tabular}{|c|c|c|}
\hline
   Signature of $(c_1,c_2,c_3)$
& Simply connected Lie group
& Property \\
  \hline
  $(+,+,+)$
& $\mathrm{SU}(2)$
& compact and simple \\
$(+,+,-)$ & $\widetilde{\mathrm{SL}}(2,\mathbb{R})$ & non-compact
and simple
\\
$(+,+,0)$ & $\widetilde{\mathrm{E}}(2)$ & solvable
\\
$(+,-,0)$ & $\mathrm{E}(1,1)$ & solvable
\\
$(+,0,0)$ & Heisenberg group  & nilpotent
\\
$(0,0,0)$ & $(\mathbb{R}^3,+)$ & Abelian
\\
\hline
\end{tabular}
\end{center}
It is easy to see that any left-invariant metric on $(\R^3,+)$ gives
rise to an isometry group of dimension 6, and in \cite{Pa} it was
proven that any left-invariant metric on the Heisenberg group gives
rise to an isometry group of dimension 4. Hence, for our purpose, we
may exclude these. The following examples cover the other cases.

\begin{example}[\textbf{The special unitary group $\mathrm{SU}(2)$}]
The group $\mathrm{SU}(2)$ is diffeomorphic to the sphere $\S^3(1)$
since $$\mathrm{SU}(2)=\left\{\left(\begin{array}{cc}
\alpha&\beta\\-\overline{\beta} &
\overline{\alpha}\end{array}\right) \ \ \biggr \vert \ \
|\alpha|^2+|\beta|^2=1 \right \}.$$ The Lie algebra of this group is
explicitly given by
$$\mathfrak{su}(2)=\left\{\left.\left(\begin{array}{cc}iu & -v+iw\\ v+iw & -iu\end{array}\right)\ \right|\ u,v,w\in\R\right\}.$$

To construct an orthonormal basis $\{e_1,e_2,e_3\}$ as mentioned
above, we proceed as follows. We take the following quaternionic
basis $\{\textbf{i},\textbf{j},\textbf{k}\}$ of $\mathfrak{su}(2)$:
$$
\textbf{i}= \left (
\begin{array}{cc}
0 & i \\
i & 0
\end{array}
\right ), \ \ \textbf{j}= \left (
\begin{array}{cc}
0 & -1 \\
1 & 0
\end{array}
\right ), \ \ \textbf{k}= \left (
\begin{array}{cc}
i & 0 \\
0 & -i
\end{array}
\right ),
$$
and we denote the left-translated vector fields of
$\{\textbf{i},\textbf{j},\textbf{k}\}$ by $\{E_1,E_2,E_3\}$. Choose
strictly positive real constants $\lambda_1$, $\lambda_2$,
$\lambda_3$ and define
$$
e_{1}=\frac{1}{\lambda_{2}\lambda_{3}} E_1,\ \
e_{2}=\frac{1}{\lambda_{3}\lambda_{1}} E_2,\ \
e_{3}=\frac{1}{\lambda_{1}\lambda_{2}} E_3.
$$
Then $ [e_1,e_2]=c_{3}e_{3}$, $[e_2,e_3]=c_{1}e_{1}$ and
$[e_3,e_1]=c_{2}e_{2}$, with $ c_{i}=2/\lambda_{i}^{2}$. The
left-invariant metric $g(c_1,c_2,c_3)$, defined by the condition
that $\{e_1,e_2,e_3\}$ is an orthonormal frame, is
$$
g(c_1,c_2,c_3)=4\left(\frac{1}{c_2c_3}\omega_{1}^{2}+
\frac{1}{c_1c_3}\omega_{2}^{2}+
\frac{1}{c_1c_2}\omega_{3}^{2}\right),
$$
where $\{\omega_1,\omega_2,\omega_3\}$ is the dual coframe field of
$\{E_1,E_2,E_3\}$.

\begin{proposition}[\cite{Pa}]\label{prop3}
Any left-invariant metric on $\mathrm{SU}(2)$ is isometric to one of
the metrics $g(c_1,c_2,c_3)$, with $c_1,c_2,c_3\geq 0$. Moreover,
the dimension of the isometry group is $\geq 4$ if and only if at
least two of the parameters $c_i$ coincide.
\end{proposition}
In particular, if $c_1=c_2=c_3=c>0$, then the space is of constant
curvature $c^2/4$.
\end{example}

\begin{example}[\textbf{The real special linear group $\mathrm{SL}(2,\R)$}]
The group $\mathrm{SL}(2,\R)$ is defined as the following subgroup
of $\mathrm{GL}(2,\R)$:
$$\mathrm{SL}(2,\R)=\left\{\left.\left(\begin{array}{cc}a&b\\c&d\end{array}\right)\ \right|\ ad-bc=1\right\}.$$
First note that this group is isomorphic to the following subgroup
of $\mathrm{GL}(2,\C)$:
$$\mathrm{SU}(1,1)=\left\{\left.\left(\begin{array}{cc}\alpha&\beta\\\overline{\beta}&\overline{\alpha}\end{array}\right)\
\right|\ |\alpha|^2-|\beta|^2=1\right\},$$ via the isomorphism
$$\mathrm{SL}(2,\R)\rightarrow
\mathrm{SU}(1,1):\left(\begin{array}{cc}a&b\\c&d\end{array}\right)\mapsto\frac{1}{2}\left(\begin{array}{cc}i&1\\1&i\end{array}\right)
\left(\begin{array}{cc}a&b\\c&d\end{array}\right)\left(\begin{array}{cc}-i&1\\1&-i\end{array}\right).$$
The Lie algebra of $\mathrm{SU}(1,1)$ is explicitly given by
$$\mathfrak{su}(1,1)=\left\{\left.\left(\begin{array}{cc}iu & v-iw\\ v+iw & -iu\end{array}\right)\ \right|\ u,v,w\in\R\right\}.$$

We take the following split-quaternionic basis of the Lie algebra
$\mathfrak{su}(1,1)$:
$$
\textbf{i}= \left (
\begin{array}{cc}
i & 0 \\
0 & -i
\end{array}
\right ), \ \ \textbf{j}^{\prime}= \left (
\begin{array}{cc}
0 & -i \\
i & 0
\end{array}
\right ), \ \ \textbf{k}^{\prime}= \left (
\begin{array}{cc}
0 & 1 \\
1 & 0
\end{array}
\right ).
$$
Denote the left-translated vector fields of $\{\textbf{j}^{\prime},
\textbf{k}^{\prime}, \textbf{i}\}$ by $\{E_1,E_2,E_3\}$. Choose
strictly positive real constants $\lambda_1$, $\lambda_2$,
$\lambda_3$ and define
$$
e_{1}=\frac{1}{\lambda_{2}\lambda_{3}} E_{1},\ \
e_{2}=\frac{1}{\lambda_{3}\lambda_{1}} E_{2},\ \
e_{3}=\frac{1}{\lambda_{1}\lambda_{2}} E_{3}.
$$
Then $ [e_1,e_2]=c_{3}e_{3}$, $[e_2,e_3]=c_{1}e_{1}$ and
$[e_3,e_1]=c_{2}e_{2}$, with $c_{1}=2/\lambda_{1}^{2}$,
$c_{2}=2/\lambda_{2}^{2}$ and $c_{3}=-2/\lambda_{3}^{2}$. The
left-invariant metric $g(c_1,c_2,c_3)$, defined by the condition
that $\{e_1,e_2,e_3\}$ is an orthonormal basis, is
$$
g(c_1,c_2,c_3)=4\left( -\frac{1}{c_{2}c_{3}}\omega_{1}^{2}
-\frac{1}{c_{3}c_{1}}\omega_{2}^{2}+
\frac{1}{c_{1}c_{2}}\omega_{3}^{2} \right),
$$
where $\{\omega_1,\omega_2,\omega_3\}$ is the dual coframe field of
$\{E_1,E_2,E_3\}$.

\begin{proposition}[\cite{Pa}] $\label{P:3.2.7}$
Any left-invariant metric on
$\mathrm{SL}(2,\R)\cong\mathrm{SU}(1,1)$ is isometric to one of the
metrics $g(c_1,c_2,c_3)$ with $c_1\geq c_2>0>c_3$. Moreover, this
metric gives rise to an isometry group of dimension $4$ if and only
if $c_1=c_2$.
\end{proposition}
\end{example}

\begin{example}[\textbf{The Minkowski motion group $\mathrm{E}(1,1)$}]
Let $\mathrm{E}(1,1)$ be the group of orientation preserving
isometries of the Minkowski plane:
$$
\mathrm{E}(1,1)=\left\{ \left(
\begin{array}{ccc}
e^{z} & 0 & x\\
0 & e^{-z} & y\\
0 & 0 & 1
\end{array}
\right) \ \Biggr \vert \ x,y,z \in \mathbb{R}\ \right\}.
$$

Consider the following left-invariant frame on $\mathrm{E}(1,1)$:
\begin{equation}\label{e'sincoordE(1,1)}e_1=\frac{1}{\lambda_1\sqrt{2}}(-e^z\partial_x+e^{-z}\partial_y),\qquad
e_2=\frac{1}{\lambda_2\sqrt{2}}(e^z\partial_x+e^{-z}\partial_y),
\qquad e_3=\frac{1}{\lambda_3}\partial_z,\end{equation} where
$\lambda_1$, $\lambda_2$ and $\lambda_3$ are strictly positive
constants. Remark that $\{e_1,e_2,e_3\}$ satisfies the commutation
relations $[e_1,e_2]=0$, $[e_2,e_3]=c_1e_1$ and $[e_3,e_1]=c_2e_2$,
with $c_1=\lambda_1/(\lambda_2\lambda_3)>0$ and
$c_2=-\lambda_2/(\lambda_1\lambda_3)<0$. We equip $\mathrm{E}(1,1)$
with a left-invariant Riemannian metric such that $\{e_1,e_2,e_3\}$
is orthonormal. The resulting Riemannian metric is
$$
g_{(\lambda_1,\lambda_2,\lambda_3)}
=\frac{\lambda_1^2}{2}(-e^{-z}dx+e^zdy)^{2}
+\frac{\lambda_2^2}{2}(e^{-z}dx+e^zdy)^{2} +\lambda_{3}^{2}dz^2.
$$
This metric is $4$-symmetric (cfr. \cite{K}) if and only if
$\lambda_1=\lambda_2$. Moreover, we have the following:
\begin{proposition}[\cite{Pa}]\label{P:3.2.3}
Any left invariant metric on $\mathrm{E}(1,1)$ is isometric to one
of the metrics $g_{(\lambda_1,\lambda_2,\lambda_3)}$ with
$\lambda_{1}\geq \lambda_{2}>0$ and
$\lambda_{3}=1/(\lambda_{1}\lambda_{2})$.
\end{proposition}
\noindent For simplicity of notation, we put
$g(\lambda_1,\lambda_2)=
g_{(\lambda_1,\lambda_2,1/(\lambda_1\lambda_2))}$.

Remark that for $\lambda_1=\lambda_2$ we have
\begin{eqnarray*}
\langle\Rt(X,Y)Z,W\rangle &=& \lambda_1^4(\langle X,W\rangle\langle Y,Z\rangle-\langle X,Z\rangle\langle Y,W\rangle\\
& & + 2\langle X,Z\rangle\langle Y,e_3\rangle\langle W,e_3\rangle + 2\langle Y,W\rangle\langle X,e_3\rangle\langle Z,e_3\rangle\\
& & - 2\langle X,W\rangle\langle Y,e_3\rangle\langle Z,e_3\rangle -
2\langle Y,Z\rangle\langle X,e_3\rangle\langle W,e_3\rangle).
\end{eqnarray*}
The Riemannian homogeneous manifold $\Sol=(\mathrm{E}(1,1),g(1,1))$
is the model space of solve-geometry in the sense of Thurston. Thus,
we have obtained the fact that $\Sol$ has a natural 2-parametric
deformation family $ \{(\mathrm{E}(1,1),g(\lambda_1,\lambda_2)) \
\vert \ \lambda_{1}\geq \lambda_{2}>0 \}. $ Note that this
deformation preserves the unimodularity property, because all these
spaces have common underlying Lie group $\mathrm{E}(1,1)$.
\end{example}

\begin{example}[\textbf{The universal covering of the Euclidean motion group $\widetilde{\mathrm{E}}(2)$}]
The group $\mathrm{E}(2)$ of orientation-preserving rigid motions of
Euclidean plane is given explictly by the following matrix group:
$$
\mathrm{E}(2)= \left\{ \left(
\begin{array}{ccc}
\cos \theta & -\sin \theta & x\\
\sin \theta & \cos \theta & y\\
0 & 0 & 1
\end{array}
\right) \ \Biggr \vert \ x,y \in \mathbb{R},\ \theta\in\S^1
\right\}.
$$
Let $\widetilde{\mathrm{E}}(2)$ denote the universal covering group
of $\mathrm{E}(2)$. Then $\widetilde{\mathrm{E}}(2)$ is isomorphic
to $\mathbb{R}^3$ with group operation
$$
(x,y,z)\ast (\overline{x},\overline{y},\overline{z})=
(x+\overline{x}\cos z-\overline{y}\sin z\ , y+\overline{x}\sin z
+\overline{y}\cos z\ , z+\overline{z}).
$$

Take strictly positive constants $\lambda_1,\lambda_2$ and
$\lambda_3$ and a left-invariant frame
$$
e_1= \frac{1}{\lambda_2} \left( -\sin z\,\partial_x+ \cos z\,
\partial_y \right), \quad e_2=\frac{1}{\lambda_3}
\partial_z, \quad e_3= \frac{1}{\lambda_1} \left
(\cos z\,\partial_x+ \sin z\,\partial_y \right).
$$
Then this frame satisfies the commutation relations $
[e_1,e_2]=c_{3}e_3$, $[e_2,e_3]=c_{1}e_1$ and $[e_3,e_1]=0$ with
$c_3=\lambda_1/(\lambda_2\lambda_3)>0$ and
$c_1=\lambda_2/(\lambda_1\lambda_3)>0$. The left-invariant
Riemannian metric determined by the condition that $\{e_1,e_2,e_3\}$
is orthonormal, is given by
$$g_{(\lambda_1,\lambda_2,\lambda_3)}=\lambda_1^2(\cos z\,dx+\sin
z\,dy)^2+\lambda_2^2(-\sin z\,dx+\cos z\,dy)^2+\lambda_3^2\,dz^2.$$
Also in this case, we have:
\begin{proposition}[\cite{Pa}]\label{P:3.2.4}
Any left-invariant metric on $\widetilde{\mathrm{E}}(2)$ is
isometric to one of the metrics
$g_{(\lambda_1,\lambda_2,\lambda_3)}$ with $\lambda_{1} >
\lambda_{2}>0$ and $\lambda_{3}=1/(\lambda_{1}\lambda_{2})$, or
$\lambda_1=\lambda_2=\lambda_3=1$. Clearly,
$\widetilde{\mathrm{E}}(2)$ with metric $g_{(1,1,1)}$ is isometric
to Euclidean three-space $\mathbb{E}^{3}$.
\end{proposition}
\noindent For simplicity of notation, we put
$g(\lambda_1,\lambda_2)=
g_{(\lambda_1,\lambda_2,1/(\lambda_1\lambda_2))}$.
\end{example}

\subsection{Non-unimodular Lie groups}
Let $G$ be a non-unimodular 3-dimensional Lie group with a
left-invariant metric. Then the \textit{unimodular kernel}
$\mathfrak{u}$ of the Lie algebra $\mathfrak{g}$ of $G$ is defined
by
$$
\mathfrak{u}=\{X \in \mathfrak{g} \ \vert \ \mathrm{tr}
(\mathrm{ad}(X))=0\}.
$$
Here $\mathrm{ad}:\mathfrak{g}\to \mathrm{End}(\mathfrak{g})$ is a
homomorphism defined by $\mathrm{ad}(X)Y=[X,Y]$. One can see that
$\mathfrak{u}$ is an ideal of $\mathfrak{g}$ which contains the
ideal $[\mathfrak{g},\mathfrak{g}]$.

It is proven in \cite{Milnor} that we can take an orthonormal basis
$\{e_1,e_2,e_3\}$ of $\mathfrak{g}$ such that
$$[e_1,e_2]=a e_{2}+b e_{3}, \quad [e_2,e_3]=0, \quad [e_1,e_3]=c e_2+d
e_3,$$ with $a+d\not=0$ and $ac+bd=0$. It is crucial that $G$ is
non-unimodular, since $e_1$ is perpendicular to $\mathfrak{u}$. In
the same article, it is remarked that after a suitable homothetic
change of the metric, we may assume that $a+d=2$. Then the constants
$a$, $b$, $c$ and $d$ are represented as
$$
a=1+\xi,\ b=(1+\xi)\eta,\ c=-(1-\xi)\eta,\ d=1-\xi,
$$
with $\xi,\eta\geq 0$. From now on, we work under this
normalization. We refer to the constants $(\xi,\eta)$ as the
\emph{structure constants} of the non-unimodular Lie group.
\begin{proposition} \label{P:3.2.8} Let $G$ be a 3-dimensional non-unimodular Lie
group with left-invariant metric $\langle\cdot,\cdot\rangle$ and use
the notations introduced above. The Levi Civita connection
$\widetilde{\nabla}$ of $G$ is given by
$$
\begin{array}{lll}
\widetilde{\nabla}_{e_1}e_{1}=0, & \widetilde{\nabla}_{e_1}e_{2}=\eta e_{3}, & \widetilde{\nabla}_{e_1}e_{3}=-\eta e_{2},\\
\widetilde{\nabla}_{e_2}e_{1}=-(1+\xi)e_{2}-\xi\eta e_{3}, &
\widetilde{\nabla}_{e_2}e_{2}=(1+\xi) e_1, &
\widetilde{\nabla}_{e_2}e_{3}=\xi\eta e_{1},\\
\widetilde{\nabla}_{e_3}e_{1}=-\xi\eta e_{2} -(1-\xi)e_{3}, &
\widetilde{\nabla}_{e_3}e_{2}=\xi\eta e_{1}, &
\widetilde{\nabla}_{e_3}e_{3}=(1-\xi)e_1.
\end{array}
$$
The frame $\{e_1,e_2,e_3\}$ diagonalizes the Ricci tensor and the
Riemann-Christoffel curvature tensor $\widetilde{R}$ is determined
by the following sectional curvatures:
\begin{eqnarray*}
&& K_{12} = \langle \widetilde{R}(e_1,e_2)e_2,e_1
\rangle=-(\xi\eta^{2}+(1+\xi)^{2}+\xi\eta^{2}(1+\xi)),\\
&& K_{23} = \langle \widetilde{R}(e_2,e_3)e_3,e_2
\rangle=\xi^{2}(1+\eta^{2})-1,\\
&& K_{13} = \langle \widetilde{R}(e_1,e_3)e_3,e_1
\rangle=\xi\eta^{2}-(1-\xi)^{2}+\xi\eta^{2}(1-\xi).
\end{eqnarray*}
\end{proposition}

Remark that for $\xi=0$, the group $G$ has constant sectional
curvature $-1$. If we assume $G$ to be simply connected, this
implies that $G$ is diffeomorphic to $\H^3(-1)$. On the other hand,
if $\xi=1$, then $G$ is locally isometric to a space with
4-dimensional isometry group. Hence, in the non-unimodular case, we
may restrict to groups with $\xi\notin\{0,1\}.$

\section{Classification results}

\subsection{In $\mathrm{SU}(2)$}
The following theorem characterizes CMC surfaces with vertically
harmonic Gauss map in $\mathrm{SU}(2)$ with 3-dimensional isometry
group.
\begin{theorem}\label{Theo5.1}
Consider the Riemannian manifold $(\mathrm{SU}(2),g(c_1,c_2,c_3))$,
with a 3-dimensional isometry group. Then we may assume that
$c_1>c_2>c_3>0$. If $f:M^2\rightarrow \mathrm{SU}(2)$ is a CMC
surface with vertically harmonic Gauss map, then $c_1=c_2+c_3$, i.e.
$K_{12}=K_{13}$, and $e_1$ is tangent to $f(M^2)$. Moreover, the
surface is minimal and the Gauss map is harmonic.
\end{theorem}

\noindent\emph{Proof.} Let
$f:M^2\rightarrow(\mathrm{SU}(2),g(c_1,c_2,c_3))$ be a  CMC surface
with vertically harmonic Gauss map and unit normal $N=\alpha
e_1+\beta e_2+\gamma e_3$, where $\{e_1,e_2,e_3\}$ is the
orthonormal frame constructed in section \ref{sec3}. It follows from
Lemma \ref{lem3.2} and Lemma \ref{lem3.3} that there are 7 cases to
consider. Cases (i), (ii), and (iii) are impossible due to the
theorem of Frobenius and formulae \eqref{Lieunimod}. From
Proposition \ref{prop2} and the assumption $c_1>c_2>c_3>0$, it
follows that only case (iv) of Lemma \ref{lem3.3} can occur.

In this case, we have $K_{12}=K_{13}$ or equivalently
\begin{equation}\label{5.1}c_1=c_2+c_3,\end{equation}
and the unit normal on $f(M^2)$ takes the form $N=\beta e_2+\gamma
e_3$, with $\beta^2+\gamma^2=1$. Then $E_1=\mathrm{d}f^{-1}(e_1)$
and $E_2=\mathrm{d}f^{-1}(-\gamma e_2+\beta e_3)$ are an orthonormal
frame on $M^2$. A straightforward computation using Proposition
\ref{prop2} yields
$$\mathrm{d}f([E_1,E_2])=\beta\gamma(c_2-c_3)\mathrm{d}f(E_2)+(\beta E_1[\gamma]-\gamma
E_1[\beta]-\beta^2c_2-\gamma^2c_3)N.$$ Hence, the distribution
spanned by $\{e_1,-\gamma e_2+\beta e_3\}$ is integrable if and only
if
\begin{equation}\label{5.2} \beta E_1[\gamma]-\gamma
E_1[\beta]=\beta^2c_2+\gamma^2c_3
\end{equation}
and in this case, we have $[E_1,E_2]=\beta\gamma(c_2-c_3)E_1$.

We can compute the shape operator $S$ associated to $N$ by using the
definition of $S$, Proposition \ref{prop2}, \eqref{5.1} and
\eqref{5.2}, to be
$$S=\left(\begin{array}{cc}0&\beta^2c_2+\gamma^2c_3\\ \beta^2c_2+\gamma^2c_3&\gamma E_2[\beta]-\beta E_2[\gamma]\end{array}\right).$$
Hence, the surface is CMC if and only if $\gamma E_2[\beta]-\beta
E_2[\gamma]=C$ is constant. Together with $\beta^2+\gamma^2=1$, we
obtain
$$\left\{\begin{array}{ll}E_2[\beta]=C\gamma,\\E_2[\gamma]=-C\beta.\end{array}\right.$$
Similarly, from \eqref{5.2} and $\beta^2+\gamma^2=1$, we obtain
$$\left\{\begin{array}{ll}E_1[\beta]=\gamma(\beta^2c_2+\gamma^2c_3),\\E_1[\gamma]=-\beta(\beta^2c_2+\gamma^2c_3).\end{array}\right.$$
By expressing the compatibility condition for the equations for
$\beta$, i.e. $E_1[E_2[\beta]]-E_2[E_1[\beta]]=[E_1,E_2][\beta]$, we
obtain $3C\beta\gamma^2(c_3-c_2)=0$. Since the distributions spanned
by $\{e_1,e_2\}$ and $\{e_1,e_3\}$ are not integrable, we have
$\beta\neq 0$ and $\gamma\neq 0$. Hence we obtain $C=0$, or
equivalently, the immersion is minimal. Remark that the
compatibility condition for the equations for $\gamma$ is then also
satisfied.

Finally, using Lemma \ref{lem3.1}, we see that the Gauss map is
harmonic. Indeed, the vector fields $U_1=(E_1+E_2)/\sqrt{2}$ and
$U_2=(E_1-E_2)/\sqrt{2}$ diagonalize the shape operator and a
straightforward computation yields $\langle
R^{\mathrm{SU}(2)}(N,\mathrm{d}f(U_i))\mathrm{d}f(U_i),N\rangle=0$
for $i=1,2$. \hfill$\square$\medskip

\subsection{In $\mathrm{SL}(2,\R)$}
The following theorem can be proven analogously as Theorem
\ref{Theo5.1}.
\begin{theorem}\label{Theo5.2}
Consider the Riemannian manifold
$(\mathrm{SL}(2,\R),g(c_1,c_2,c_3))$, with a 3-dimensional isometry
group. Then we may assume that $c_1>c_2>0>c_3$. If $f:M^2\rightarrow
\mathrm{SL}(2,\R)$ is a CMC surface with vertically harmonic Gauss
map, then $c_2=c_1+c_3$, i.e. $K_{12}=K_{23}$ and $e_2$ is tangent
to $f(M^2)$. Moreover, the surface is minimal and the Gauss map is
harmonic.
\end{theorem}

\subsection{In $\mathrm{E}(1,1)$}
Since formulae \eqref{e'sincoordE(1,1)} give the orthonormal frame
field $\{e_1,e_2,e_3\}$ on
$(\mathrm{E}(1,1),g(\lambda_1,\lambda_2))$ explicitly in terms of
the natural coordinate vector fields
$\{\partial_x,\partial_y,\partial_z\}$, we can describe the CMC
surfaces with vertically harmonic Gauss map explicitly in these
coordinates.

\begin{theorem}\label{Theo5.3}
Consider the Riemannian manifold
$(\mathrm{E}(1,1),g(\lambda_1,\lambda_2))$ and let
$$f:U\subseteq\R^2\rightarrow \mathrm{E}(1,1):(u,v)\mapsto\left(\begin{array}{ccc}e^{f_3(u,v)} & 0 & f_1(u,v)\\0 & e^{-f_3(u,v)} & f_2(u,v)\\0 & 0 & 1\end{array}\right)$$
be a CMC surface with vertically harmonic Gauss map. Then the Gauss
map is harmonic, the surface is minimal and there exists an open
subset $V\subseteq U$ on which $f$ is, up to reparametrization and
isometries of the ambient space, given by one of the following:
\begin{itemize}
\item[(i)] $f_1=u$, $f_2=v$ and $f_3=0$,
\item[(ii)] $f_1=u$, $f_2=-u$ and $f_3=v$.
\end{itemize}
The latter case only occurs if $\lambda_1=\lambda_2$. The surfaces
of the first type are flat. Conversely, all surfaces described above
are minimal and have harmonic Gauss map.
\end{theorem}

\noindent\emph{Proof.} Let $f:U\subseteq\R^2\rightarrow
(\mathrm{E}(1,1),g(\lambda_1,\lambda_2))$ be a CMC surface with
vertically harmonic Gauss map and suppose that $N=\alpha e_1+\beta
e_2+\gamma e_3$ is a unit normal on the surface. From Lemma
\ref{lem3.2}, Lemma \ref{lem3.3}, the theorem of Frobenius, formulae
\eqref{Lieunimod} and Proposition \ref{prop2}, it follows that there
are only two cases to consider, namely `$\alpha=\beta=0$' and
`$\gamma=0$ and $K_{13}=K_{23}$'.

If $\alpha=\beta=0$, we may assume that $N=e_3$. A straightforward
calculation shows that
$$S=-\frac{\lambda_1^2+\lambda_2^2}{2}\left(\begin{array}{cc} 0&1\\1&0 \end{array}\right)$$
with respect to the basis
$\{E_1=\mathrm{d}f^{-1}(e_1),E_2=\mathrm{d}f^{-1}(e_2)\}$. Hence the
surface is minimal, such that Lemma \ref{lem3.2} implies that it has
vertically harmonic Gauss map. Since
$\mathrm{span}\{e_1,e_2\}=\mathrm{span}\{\partial_x,\partial_y\}$,
this case corresponds, after an isometry of the ambient space, to
the first surface given in the theorem. Remark that, using the
equation of Gauss, the Gaussian curvature of the surface is given by
$K=\langle R^{E(1,1)}(e_1,e_2)e_2,e_1\rangle + \det S = 0.$ From
Lemma 1, we obtain that the Gauss map is harmonic. Indeed, the
vector fields $U_1=(E_1+E_2)/\sqrt{2}$ and $U_2=(E_1-E_2)/\sqrt{2}$
form an orthonormal frame which diagonalizes the shape operator and
$\langle
R^{E(1,1)}(N,\mathrm{d}f(U_i))\mathrm{d}f(U_i),N\rangle=(K_{13}+K_{23})/2$
for $i=1,2$.

In the second case, it follows from $K_{13}=K_{23}$ that
$\lambda_1=\lambda_2$. Moreover, the tangent plane to $M^2$ is at
every point spanned by the orthonormal vector fields
$E_1=\mathrm{d}f^{-1}(\beta e_1-\alpha e_2$) and
$E_2=\mathrm{d}f^{-1}(e_3)$. The Lie bracket of these vector
fields is given by
$$\mathrm{d}f([E_1,E_2])=-2\lambda_1^2\alpha\beta \mathrm{d}f(E_1)+(\beta E_2[\alpha]-\alpha E_2[\beta]-\lambda_1^2(\alpha^2-\beta^2))N.$$
Hence, the integrability condition for the distribution spanned by
$\{\beta e_1-\alpha e_2,e_3\}$ is
\begin{equation}\label{2}
\beta E_2[\alpha]-\alpha E_2[\beta]=\lambda_1^2(\alpha^2-\beta^2)
\end{equation}
and we have $[E_1,E_2]=-2\lambda_1^2\alpha\beta E_1$.

The shape operator associated to $N=\alpha e_1+\beta e_2$ with
respect to the basis $\{E_1,E_2\}$ is given by
$$S=\left(\begin{array}{cc}\alpha E_1[\beta]-\beta E_1[\alpha] &
\lambda_1^2(\beta^2-\alpha^2)\\
\lambda_1^2(\beta^2-\alpha^2) & 0
\end{array}\right).$$

We are interested in the case that $M^2$ is a CMC surface. Hence
assume that
\begin{equation}\label{3}
\alpha E_1[\beta]-\beta E_1[\alpha]=C
\end{equation}
for some real constant $C$. The fact that $\alpha^2+\beta^2=1$,
together with (\ref{2}) and (\ref{3}) gives the following equations:
$E_1[\alpha] = -C\beta$, $E_1[\beta] = C\alpha$, $E_2[\alpha] =
\lambda_1^2\beta(\alpha^2-\beta^2)$, $E_2[\beta] =
-\lambda_1^2\alpha(\alpha^2-\beta^2)$. If we introduce a function
$\theta$ which is locally defined on the surface by
$\alpha=\cos\theta$, $\beta=\sin\theta$, these equations reduce to
\begin{equation}\label{6}
\left\{\begin{array}l
E_1[\theta]=C,\\E_2[\theta]=-\lambda_1^2\cos(2\theta).\end{array}\right.
\end{equation}
The compatibility condition for this system is
$C\lambda_1^2\sin(2\theta)=0$, from which we conclude $C=0$.
Hence, the surface is minimal.

We can take coordinates $(u,v)$ on the surface, with
$\partial_u=E_1$ and $\partial_v=pE_1+E_2$ for a suitable function
$p:M^2\rightarrow\R$. Indeed, the condition
$[\partial_u,\partial_v]=0$ is equivalent to the equation
\begin{equation}\label{5}
\partial_u p=\lambda_1^2\sin(2\theta).
\end{equation}
Moreover, the system of equations (\ref{6}) is equivalent to
$$\left\{\begin{array}{ll}\partial_u\theta = 0,\\
\partial_v\theta =
-\lambda_1^2\cos(2\theta),\end{array}\right.$$ which can be solved
as
$$\theta=\arctan\left(\frac{e^{-2\lambda_1^2v+c}+1}{e^{-2\lambda_1^2v+c}-1}\right),$$
where $c$ is a real constant. Remark that (\ref{5}) yields
$p(u,v)=\lambda_1^2u\sin(2\theta(v))+C(v)$ for some function
$C(v)$. Since we are only interested in one coordinate system, we
may assume that $C(v)=0$ and hence
$$p(u,v)=\lambda_1^2u\sin(2\theta(v)).$$

In order to find an explicit expression for $f$, we need to
integrate the formulae
\begin{eqnarray*}
& & \left(\partial_u f_1, \partial_u f_2, \partial_u f_3\right) =
\mathrm{d}f(E_1) = \sin\theta e_1
-\cos\theta e_2,\\
& & \left(\partial_v f_1, \partial_vf_2, \partial_v f_3\right) = p
\mathrm{d}f(E_1)+\mathrm{d}f(E_2) = p\sin\theta e_1-p\cos\theta
e_2 + e_3,
\end{eqnarray*}
with
\begin{eqnarray*}
e_1 &=& \frac{1}{\lambda_1\sqrt{2}}(-e^{f_3},e^{-f_3},0),\\
e_2 &=& \frac{1}{\lambda_1\sqrt{2}}(e^{f_3},e^{-f_3},0),\\
e_3 &=& \lambda_1^2(0,0,1).
\end{eqnarray*}
A direct computation yields
\begin{eqnarray*}
f_1 &=& -\frac{ue^{\lambda_1^2v}}{\lambda_1\sqrt{2}}(\sin\theta(v)+\cos\theta(v))+a_1,\\
f_2 &=& \frac{ue^{-\lambda_1^2v}}{\lambda_1\sqrt{2}}(\sin\theta(v)-\cos\theta(v))+a_2,\\
f_3 &=& \lambda_1^2v,
\end{eqnarray*}
where $a_1$ and $a_2$ are real constants. The left translation $L_A
: \mathrm{E}(1,1) \to \mathrm{E}(1,1)$, with
$$ A = \left( \begin{array}{ccc} e^{c/2} & 0 & -a_1 \\ 0 & e^{-c/2} & -a_2 \\ 0 & 0 & 1 \end{array} \right),$$
which is of course an isometry of $(\mathrm{E}(1,1),
g(\lambda_1,\lambda_1))$, maps the image of $f$ into the surface
given by $x=-y$. This corresponds to the second case given in the
theorem.

The vector fields $U_1=(E_1+E_2)/\sqrt{2}$ and
$U_2=(E_1-E_2)/\sqrt{2}$ form an orthonormal frame which
diagonalizes the shape operator. It follows from Lemma 1 and the
observation $\langle R^{\mathrm{E}(1,1)}(N,\mathrm{d}f(U_i))$
$\mathrm{d}f(U_i),N\rangle=0$ for $i=1,2$, that the Gauss map is
harmonic.

Conversely, one can check that both surfaces given in the theorem
are minimal and satisfy the conditions of Lemma \ref{lem3.1}. Hence
they have harmonic Gauss maps. \hfill$\square$\medskip

As a corollary, we obtain that there are two families of CMC
surfaces in $\mathrm{Sol}_3$ with harmonic Gauss map. Both of them
are minimal and explicit parametrizations are given in Theorem
\ref{Theo5.3}, if we put $\lambda_1=\lambda_2=1$.

\subsection{In $\widetilde{\mathrm{E}}(2)$}
Also in this case, we can give the solutions to our classification
problem by means of an explicit parametrization. The proof is again
similar.

\begin{theorem}\label{Theo5.4} Up to isometries of the ambient space, all CMC surfaces in $(\widetilde{\mathrm{E}}(2),g(\lambda_1,\lambda_2))$
with vertically harmonic tangential Gauss map can be locally
parametrized as
$$f:U\subseteq\R^2\rightarrow\widetilde{\mathrm{E}}(2):(u,v)\mapsto\left(\begin{array}{ccc}1&0&u\\0&1&v\\0&0&1\end{array}\right).$$
All these surfaces are minimal and their Gauss map is harmonic.
\end{theorem}

\subsection{In non-unimodular Lie groups}
As mentioned at the end of section \ref{sec3}, we may restrict
ourselves to non-unimodular groups with $\xi\notin\{0,1\}$. The
following theorem can be proven again with similar methods.
\begin{theorem}
Let $G$ be a 3-dimensional non-unimodular Lie group with
left-invariant metric. We may assume that the first structure
constant satisfies $\xi\notin\{0,1\}$. If $f:M^2\rightarrow G$ is a
CMC surface with vertically harmonic Gauss map, then, up to
isometries of $G$, one of the following holds:
\begin{itemize}
\item[(i)] $f(M^2)$ is an integral surface of the distribution spanned by
$\{e_2,e_3\}$. The surface is flat and has constant mean curvature
$H=1$.
\item[(ii)] $f(M^2)$ is an integral surface of the distribution spanned by
$\{e_1,e_2\}$ or of the distribution spanned by $\{e_1,e_3\}$. These
surfaces are totally geodesic and hence minimal.
\end{itemize}
The second case only occurs if $\eta=0$.
\end{theorem}

\end{document}